\documentclass[11pt]{article}
\usepackage{amsfonts}
\usepackage{amsfonts}
\usepackage{amsfonts}
\usepackage{mathrsfs}
\usepackage{bm}
\usepackage{cite}
\usepackage[colorlinks,linkcolor=blue,citecolor=blue]{hyperref}
\usepackage{amssymb,amsmath} 
 \allowdisplaybreaks

\catcode`!=11
\let\!int\int \def\int{\displaystyle\!int}
\let\!lim\lim \def\lim{\displaystyle\!lim}
\let\!sum\sum \def\sum{\displaystyle\!sum}
\let\!sup\sup \def\sup{\displaystyle\!sup}
\let\!inf\inf \def\inf{\displaystyle\!inf}
\let\!cap\cap \def\cap{\displaystyle\!cap}
\let\!max\max \def\max{\displaystyle\!max}
\let\!min\min \def\min{\displaystyle\!min}
\let\!frac\frac \def\frac{\displaystyle\!frac}
\catcode`!=12

\let\oldsection\section
\renewcommand\section{\setcounter{equation}{0}\oldsection}

\allowdisplaybreaks
\def\pf{\it{Proof.}\rm\quad}

\def\N{\mathbb{N}}\def\Z{\mathbb{Z}}

\newtheorem{thm}{Theorem}[section]
\newtheorem{lem}[thm]{Lemma}
\newtheorem{cor}[thm]{Corollary}

\begin{document}
\title {\bf Some infinite series involving hyperbolic functions}
\author{
{Ce Xu\thanks{Corresponding author. Email: xuce1242063253@163.com}}\\[1mm]
\small School of Mathematical Sciences, Xiamen University\\
\small Xiamen
361005, P.R. China}

\date{}
\maketitle \noindent{\bf Abstract} This paper develops an approach to the evaluation of infinite series involving hyperbolic functions. By using the approach, we give explicit formulae for several classes of series of hyperbolic functions in terms of Riemann zeta values. Moreover, we also establish many relations involving two or more series of hyperbolic functions. Furthermore, we obtain the Ramanujan's formula for $\zeta(2n+1)$ and find another similar formulae. The approach is based on simple contour integral representations and residue computations. Some interesting (known or new) consequences and illustrative examples are considered.
\\[2mm]
\\[2mm]
\tableofcontents
\section{Introduction}
Extensive use is made throughout the paper of the Bernoulli and Euler numbers $B_n$ and $E_n$ that are defined in terms of the Bernoulli and Euler polynomials of order $n$, $B_n(x)$ and $E_n(x)$, respectively. These polynomials are defined by the generating functions \cite{GR2007}
\begin{align*}
\frac{{t{e^{xt}}}}
{{{e^t} - 1}} = \sum\limits_{n = 0}^\infty  {{B_n}\left( x \right)\frac{{{t^n}}}
{{n!}}}\quad {\rm for}\quad \left| t \right| < 2\pi
\end{align*}
and
\begin{align*}
\frac{{2{e^{xt}}}}
{{{e^t} + 1}} = \sum\limits_{n = 0}^\infty  {{E_n}\left( x \right)\frac{{{t^n}}}
{{n!}}}\quad {\rm for}\quad \left| t \right| < \pi.
\end{align*}
The Bernoulli numbers are always denoted by $B_n$ and are defined by the relation
\[B_n=B_n(0)\quad {\rm for}\quad n=0,1,2,\ldots, \]
when
\[B_0=1,\ B_1=-\frac 1{2}, B_2=\frac 1{6}, B_4=-\frac 1{30}, B_6=\frac 1{42}, \ldots.\]
The Euler numbers $E_n$ are defined by setting
\[E_n=2^nE_n\left(\frac 1{2}\right)\quad {\rm for} \quad n=0,1,2,\ldots\]
The $E_n$ are all integral, and $E_0=1,\ E_2=-1,\ E_4=5,\ E_6=-61,\ldots$.

The subject of this paper is hyperbolic functions. The basic hyperbolic functions are the hyperbolic sine "sinh" and the hyperbolic cosine "cosh"  from which are derived the hyperbolic tangent "tanh", hyperbolic cosecant "csch" or "cosech" , hyperbolic secant "sech", and hyperbolic cotangent "coth", corresponding to the derived trigonometric functions. They may be defined in terms of the exponential function:
\begin{align*}
 & \sinh x = \frac{{{e^x} - {e^{ - x}}}}
{2}, \hfill \\
  &\cosh x = \frac{{{e^x} + {e^{ - x}}}}
{2}, \hfill \\
  &\tanh x = \frac{{\sinh x}}
{{\cosh x}} = \frac{{{e^x} - {e^{ - x}}}}
{{{e^x} + {e^{ - x}}}}, \hfill \\
  &\coth x = \frac{{\cosh x}}
{{\sinh x}} = \frac{{{e^x} + {e^{ - x}}}}
{{{e^x} - {e^{ - x}}}},\quad x\neq 0, \hfill \\
  &{\rm sech}x = \frac{1}
{{\cosh x}} = \frac{2}
{{{e^x} + {e^{ - x}}}}, \hfill \\
  &{\rm csch}x = \frac{1}
{{\sinh x}} = \frac{2}
{{{e^x} - {e^{ - x}}}},\quad x\neq 0,.
\end{align*}
Infinite series involving the hyperbolic functions have attracted the attention of many authors. Berndt's papers \cite{B1977,B1978,B2004,B2016} and paper \cite{BS2017} with A. Straub, and his books \cite{B1985,B1989} and his fourth book with Andrews \cite{AB2013} contain many such results well as numerous references.

The purposes of this paper are to establish some relations between Riemann zeta function and infinite series involving hyperbolic functions. Furthermore, we evaluate some partial fraction expansions of several quadratic products involving hyperbolic and trigonometric function, and give explicit formulae. Moreover, we obtain the Ramanujan's formula for $\zeta(2n+1)$ and a similar identity. Finally, we evaluate two series involving harmonic number and hyperbolic function, and establish two general formulae. The results which we present here can be seen as an extension of Berndt's work. As a result, we attain many expressions for the combined infinite series involving hyperbolic functions by the Riemann zeta values, which contain most of the known results and some new relations. For example, we obtain
\[\sum\limits_{n = 1}^\infty  {\frac{{{\pi ^2}{{\left( { - 1} \right)}^n}}}
{{{{\sinh }^2}\left( {\pi n} \right)}}}  + \sum\limits_{n = 1}^\infty  {\frac{{{\pi ^2}\cosh \left( {\pi n} \right)}}
{{{{\sinh }^2}\left( {\pi n} \right)}}}  + \frac{1}
{2}\zeta \left( 2 \right) = 0.\]

In this paper we will use the contour integration and residue theorem to prove our main results. The contour integration is a classical technique for evaluating infinite sums by reducing them to a finite number of residue computations. This summation mechanism is formalized by a lemma that goes back to Cauchy and is nicely developed throughout \cite{FS1998}.
 
The Riemann zeta function $\zeta(s)$ is defined by
$$\zeta(s):=\sum\limits_{n = 1}^\infty {\frac {1}{n^{s}}},\Re(s)>1.$$
When $s=2k\ (k\in \N)$ is an even, Euler gave the famous formula
\[\zeta \left( {2k} \right) = \frac{{{{\left( { - 1} \right)}^{k - 1}}{B_{2k}}}}
{{2\left( {2k} \right)!}}{\left( {2\pi } \right)^{2k}}.\]
The formula not only provides an elegant formula for evaluating $\zeta(2k)$, but it also tells us of the arithmetical nature of $\zeta(2k)$. In contrast, we know very little about $\zeta(2k+1)$. One of the major achievements in number theory in the past half-century are R. Ap$\acute{\rm e}$ry's proof that $\zeta(3)$ is irrational \cite{Ap1981}, but for $k\geq 2$, the arithmetical nature of $\zeta(2k+1)$ remains open. Further results for $k\geq 2$, please see Zudilin's paper \cite{Z2001,2Z2001}.

The following lemma will be useful in the development of the main theorems.
\begin{lem}(\cite{FS1998})\label{lem1}
Let $\xi \left( s \right)$ be a kernel function and let $r(s)$ be a function which is $O(s^{-2})$ at infinity. Then
\begin{align}\label{1.1}
\sum\limits_{\alpha  \in O} {{\mathop{\rm Re}\nolimits} s{{\left( {r\left( s \right)\xi \left( s \right)} \right)}_{s = \alpha }}}  + \sum\limits_{\beta  \in S} {{\mathop{\rm Re}\nolimits} s{{\left( {r\left( s \right)\xi \left( s \right)} \right)}_{s = \beta }}}  = 0.
\end{align}
where $S$ is the set of poles of $r(s)$ and $O$ is the set of poles of $\xi \left( s \right)$ that are not poles $r(s)$ . Here ${\mathop{\rm Re}\nolimits} s{\left( {r\left( s \right)} \right)_{s = \alpha }} $ denotes the residue of $r(s)$ at $s= \alpha$. The kernel function $\xi \left( s \right)$ is meromorphic in the whole complex plane and satisfies $\xi \left( s \right)=o(s)$ over an infinite collection of circles $\left| z \right| = {\rho _k}$ with ${\rho _k} \to \infty . $
\end{lem}
\section{Five general Theorems}
In \cite{FS1998}, P. Flajolet and B. Salvy gave the following asymptotic formulae
\begin{align}\label{2.1}
\pi \cot \left( {\pi z} \right)\mathop  = \limits^{z \to n} \frac{1}{{z - n}} - 2\sum\limits_{k = 1}^\infty  {\zeta \left( {2k} \right){{\left( {z - n} \right)}^{2k - 1}}} ,
\end{align}
\begin{align}\label{2.2}
\frac{\pi }
{{\sin \left( {\pi z} \right)}}\mathop  = \limits^{z \to n} {\left( { - 1} \right)^n}\left( {\frac{1}
{{z - n}} + 2\sum\limits_{k = 1}^\infty  {\bar \zeta \left( {2k} \right){{\left( {z - n} \right)}^{2k - 1}}} } \right),
\end{align}
where $n$ is integer (namely $n=0, \pm 1, \pm 2, \pm 3, \cdots $), and ${\bar \zeta} (s)$ denotes the alternating Riemann zeta function defined by
\[\bar \zeta \left( s \right) := \sum\limits_{n = 1}^\infty  {\frac{{{{\left( { - 1} \right)}^{n - 1}}}}{{{n^s}}}}=(1-2^{1-s})\zeta(s) ,\;{\mathop{\Re}\nolimits} \left( s \right) \ge 1.\]
Then using the basic relation $\sin x=\cos\left(\frac {\pi}{2}-x\right)$ and the following relations between trigonometric functions and hyperbolic functions
\begin{align*}
&\sin x =  - i\sinh \left( {ix} \right),\sinh x =  - i\sin \left( {ix} \right),\cos x = \cosh \left( {ix} \right),\cosh x = \cos \left( {ix} \right), \hfill \\
&\tan x =  - i\tanh \left( {ix} \right),\tanh x =  - i\tan \left( {ix} \right),\cot x = i\coth \left( {ix} \right),\coth x = i\cot \left( {ix} \right),
\end{align*}
by a simple calculation, the following identities are easily derived
\begin{align}\label{2.3}
\pi \coth \left( {\pi z} \right)\mathop  = \limits^{z \to ni} \frac{1}
{{z - ni}} - 2\sum\limits_{k = 1}^\infty  {{{\left( { - 1} \right)}^k}\zeta \left( {2k} \right){{\left( {z - ni} \right)}^{2k - 1}}} ,
\end{align}
\begin{align}\label{2.4}
\frac{\pi }
{{\sinh \left( {\pi z} \right)}}\mathop  = \limits^{z \to ni} {\left( { - 1} \right)^n}\left( {\frac{1}
{{z - ni}} + 2\sum\limits_{k = 1}^\infty  {{{\left( { - 1} \right)}^k}\bar \zeta \left( {2k} \right){{\left( {z - ni} \right)}^{2k - 1}}} } \right),
\end{align}
\begin{align}\label{2.5}
 \pi \tan \left( {\pi z} \right)\mathop  = \limits^{z \to \frac{{2n - 1}}
{2}}  - \frac{1}
{{z - \frac{{2n - 1}}
{2}}} + 2\sum\limits_{k = 1}^\infty  {\zeta \left( {2k} \right){{\left( {z - \frac{{2n - 1}}
{2}} \right)}^{2k - 1}}} ,
\end{align}
\begin{align}\label{2.6}
\pi \tanh \left( {\pi z} \right)\mathop  = \limits^{z \to \frac{{2n - 1}}
{2}i} \frac{1}
{{z - \frac{{2n - 1}}
{2}i}} - 2\sum\limits_{k = 1}^\infty  {{{\left( { - 1} \right)}^k}\zeta \left( {2k} \right){{\left( {z - \frac{{2n - 1}}
{2}i} \right)}^{2k - 1}}} ,
\end{align}
\begin{align}\label{2.7}
\frac{\pi }
{{\cos \left( {\pi z} \right)}}\mathop  = \limits^{z \to \frac{{2n - 1}}
{2}} {\left( { - 1} \right)^n}\left\{ {\frac{1}
{{z - \frac{{2n - 1}}
{2}}} + 2\sum\limits_{k = 1}^\infty  {\bar \zeta \left( {2k} \right){{\left( {z - \frac{{2n - 1}}
{2}} \right)}^{2k - 1}}} } \right\},
\end{align}
\begin{align}\label{2.8}
&   \frac{\pi }
{{\cosh \left( {\pi z} \right)}}\mathop  = \limits^{z \to \frac{{2n - 1}}
{2}i} {\left( { - 1} \right)^n}i\left\{ {\frac{1}
{{z - \frac{{2n - 1}}
{2}i}} + 2\sum\limits_{k = 1}^\infty  {{{\left( { - 1} \right)}^k}\bar \zeta \left( {2k} \right){{\left( {z - \frac{{2n - 1}}
{2}i} \right)}^{2k - 1}}} } \right\}.
\end{align}
Next, we will use the formulae (\ref{2.1})-(\ref{2.8}) and Lemma \ref{lem1} to prove our main results.

\subsection{Main theorems and Corollaries}
Now we state our main results. The main results of this section are the following five theorems.
\begin{thm}\label{thm2.1} For any positive integer $k>1$ and reals $x,y$ with $x,y\neq 0$, then
\begin{align}\label{2.9}
&\sum\limits_{n = 1}^\infty  {\frac{{\pi \coth \left( {\frac{{\pi ny}}
{x}} \right)}}
{{{n^{2k - 1}}}}} {x^{2k - 1}}y + {\left( { - 1} \right)^k}\sum\limits_{n = 1}^\infty  {\frac{{\pi \coth \left( {\frac{{\pi nx}}
{y}} \right)}}
{{{n^{2k - 1}}}}} {y^{2k - 1}}x\nonumber\\
& =  \left( {{x^{2k}} + {{\left( { - 1} \right)}^k}{y^{2k}}} \right)\zeta \left( {2k} \right) - 2\sum\limits_{\scriptstyle {k_1} + {k_2} = k,\hfill \atop \scriptstyle {k_1}, {k_2} \geqslant 1\hfill} {{{\left( { - 1} \right)}^{{k_2}}}{x^{2{k_1}}}{y^{2{k_2}}}\zeta \left( {2{k_1}} \right)\zeta \left( {2{k_2}} \right)}.
\end{align}
\end{thm}
\begin{thm}\label{thm2.2} For any positive integer $k$ and reals $x,y$ with $x,y\neq 0$, then
\begin{align}\label{2.10}
&\sum\limits_{n = 1}^\infty  {\frac{{\pi {{\left( { - 1} \right)}^{n - 1}}}}
{{{n^{2k - 1}}\sinh \left( {\frac{{\pi ny}}
{x}} \right)}}} {x^{2k - 1}}y + {\left( { - 1} \right)^k}\sum\limits_{n = 1}^\infty  {\frac{{\pi {{\left( { - 1} \right)}^{n - 1}}}}
{{{n^{2k - 1}}\sinh \left( {\frac{{\pi nx}}
{y}} \right)}}} {y^{2k - 1}}x\nonumber\\
& = \left( {{x^{2k}} + {{\left( { - 1} \right)}^k}{y^{2k}}} \right)\bar \zeta \left( {2k} \right) + 2\sum\limits_{\scriptstyle {k_1} + {k_2} = k,\hfill \atop \scriptstyle {k_1}, {k_2} \geqslant 1\hfill} {{{\left( { - 1} \right)}^{{k_2}}}{x^{2{k_1}}}{y^{2{k_2}}}\bar \zeta \left( {2{k_1}} \right)\bar \zeta \left( {2{k_2}} \right)}.
\end{align}
\end{thm}
\begin{thm}\label{thm2.3} For any positive integer $k$ and reals $x,y$ with $x,y\neq 0$, then
\begin{align}\label{2.11}
&\sum\limits_{n = 1}^\infty  {\frac{\pi }
{{{n^{2k - 1}}\sinh \left( {\frac{{\pi ny}}
{x}} \right)}}} {x^{2k - 1}}y + {\left( { - 1} \right)^k}\sum\limits_{n = 1}^\infty  {\frac{{\pi \coth \left( {\frac{{\pi nx}}
{y}} \right)}}
{{{n^{2k - 1}}}}{{\left( { - 1} \right)}^n}} {y^{2k - 1}}x \nonumber\\
& = {x^{2k}}\zeta \left( {2k} \right) - {\left( { - 1} \right)^k}{y^{2k}}\bar \zeta \left( {2k} \right) + 2\sum\limits_{\scriptstyle {k_1} + {k_2} = k,\hfill \atop \scriptstyle {k_1}, {k_2} \geqslant 1\hfill} {{{\left( { - 1} \right)}^{{k_2}}}{x^{2{k_1}}}{y^{2{k_2}}}\zeta \left( {2{k_1}} \right)\bar \zeta \left( {2{k_2}} \right)} .
\end{align}
\end{thm}
\begin{thm}\label{thm2.4} For any positive integer $k>1$ and reals $x,y$ with $x,y\neq 0$, then
\begin{align}\label{2.12}
&\sum\limits_{n = 1}^\infty  {\frac{{\pi \tanh \left( {\frac{{2n - 1}}
{{2x}}\pi y} \right)}}
{{{{\left( {2n - 1} \right)}^{2k - 1}}}}} {x^{2k - 1}}y + {\left( { - 1} \right)^k}\sum\limits_{n = 1}^\infty  {\frac{{\pi \tanh \left( {\frac{{2n - 1}}
{{2y}}\pi x} \right)}}
{{{{\left( {2n - 1} \right)}^{2k - 1}}}}} {y^{2k - 1}}x\nonumber\\
& =\sum\limits_{\scriptstyle {k_1} + {k_2} = k,\hfill \atop \scriptstyle {k_1}, {k_2} \geqslant 1\hfill} {\frac{{\left( {{2^{2{k_1}}} - 1} \right)\left( {{2^{2{k_2}}} - 1} \right)\left| {{B_{2{k_1}}}} \right|{B_{2{k_2}}}}}
{{\left( {2{k_1}} \right)!\left( {2{k_2}} \right)!}}{x^{2{k_1}}}{y^{2{k_2}}}{\pi ^{2k}}}.
\end{align}
\end{thm}
\begin{thm}\label{thm2.5} For any integer $k\geq 0$ and reals $x,y$ with $x,y\neq 0$, then
\begin{align}\label{2.13}
&\sum\limits_{n = 1}^\infty  {\frac{{\pi {{\left( { - 1} \right)}^{n - 1}}}}
{{{{\left( {2n - 1} \right)}^{2k + 1}}\cosh \left( {\frac{{2n - 1}}
{{2x}}\pi y} \right)}}} {x^{2k}} + {\left( { - 1} \right)^k}\sum\limits_{n = 1}^\infty  {\frac{{\pi {{\left( { - 1} \right)}^{n - 1}}}}
{{{{\left( {2n - 1} \right)}^{2k + 1}}\cosh \left( {\frac{{2n - 1}}
{{2y}}\pi x} \right)}}} {y^{2k}}\nonumber\\
&= \sum\limits_{\scriptstyle {k_1} + {k_2} = k,\hfill \atop \scriptstyle {k_1}, {k_2} \geqslant 0\hfill} {\frac{{\left| {{E_{2{k_1}}}} \right|{E_{2{k_2}}}}}
{{\left( {2{k_1}} \right)!\left( {2{k_2}} \right)!}}{x^{2{k_1}}}{y^{2{k_2}}}{{\left( {\frac{\pi }
{2}} \right)}^{2k + 2}}}.
\end{align}
\end{thm}

Setting $\alpha  = \frac{{\pi y}}{x}$ and $\beta  = \frac{{\pi \alpha }}{y}$ (namely $\alpha \beta  = {\pi ^2}$) in Theorem \ref{thm2.1}-\ref{thm2.5}, we can readily deduce the following corollaries.
\begin{cor}(\cite{N1951})\label{cor2.6} If $k>1$ is a positive integer number and $\alpha,\beta$ are real numbers such as $\alpha \beta  = {\pi ^2}$, then
\begin{align}\label{2.14}
&\alpha {\beta ^k}\sum\limits_{n = 1}^\infty  {\frac{{\coth \left( {n\alpha } \right)}}
{{{n^{2k - 1}}}}}  + {\left( { - 1} \right)^k}{\alpha ^k}\beta \sum\limits_{n = 1}^\infty  {\frac{{\coth \left( {n\beta } \right)}}
{{{n^{2k - 1}}}}}\nonumber \\
& = \left( {{\beta ^k} + \left( { - 1} \right)^k{\alpha ^k}} \right)\zeta \left( {2k} \right) - 2\sum\limits_{\scriptstyle {k_1} + {k_2} = k,\hfill \atop \scriptstyle {k_1}, {k_2} \geqslant 1\hfill}  {{{\left( { - 1} \right)}^{{k_2}}}{\alpha ^{{k_2}}}{\beta ^{{k_1}}}\zeta \left( {2{k_1}} \right)\zeta \left( {2{k_2}} \right)}.
\end{align}
\end{cor}
\begin{cor} \label{cor2.7} If $k$ is a positive integer number and $\alpha,\beta$ are real numbers such as $\alpha \beta  = {\pi ^2}$, then
\begin{align}\label{2.15}
&\alpha {\beta ^k}\sum\limits_{n = 1}^\infty  {\frac{{{{\left( { - 1} \right)}^{n - 1}}}}
{{{n^{2k - 1}}\sinh \left( {n\alpha } \right)}}}  + {\left( { - 1} \right)^k}{\alpha ^k}\beta \sum\limits_{n = 1}^\infty  {\frac{{{{\left( { - 1} \right)}^{n - 1}}}}
{{{n^{2k - 1}}\sinh \left( {n\beta } \right)}}}\nonumber \\
& = \left( {{\beta ^k} + \left( { - 1} \right)^k{\alpha ^k}} \right)\bar \zeta \left( {2k} \right) + 2\sum\limits_{\scriptstyle {k_1} + {k_2} = k,\hfill \atop \scriptstyle {k_1}, {k_2} \geqslant 1\hfill}{{{\left( { - 1} \right)}^{{k_2}}}{\alpha ^{{k_2}}}{\beta ^{{k_1}}}\bar \zeta \left( {2{k_1}} \right)\bar \zeta \left( {2{k_2}} \right)}.
\end{align}
\end{cor}
\begin{cor}\label{cor2.8} If $k$ is a positive integer number and $\alpha,\beta$ are real numbers such as $\alpha \beta  = {\pi ^2}$, then
\begin{align}\label{2.16}
&\alpha {\beta ^k}\sum\limits_{n = 1}^\infty  {\frac{1}
{{{n^{2k - 1}}\sinh \left( {n\alpha } \right)}}}  + {\left( { - 1} \right)^k}{\alpha ^k}\beta \sum\limits_{n = 1}^\infty  {\frac{{\coth \left( {n\beta } \right)}}
{{{n^{2k - 1}}}}{{\left( { - 1} \right)}^n}}\nonumber \\
&= \left( {{\beta ^k} - {{\left( { - 1} \right)}^k}\left( {1 - {2^{1 - 2k}}} \right){\alpha ^k}} \right)\zeta \left( {2k} \right) + 2\sum\limits_{\scriptstyle {k_1} + {k_2} = k,\hfill \atop \scriptstyle {k_1}, {k_2} \geqslant 1\hfill} {{{\left( { - 1} \right)}^{{k_2}}}{\alpha ^{{k_2}}}{\beta ^{{k_1}}}\zeta \left( {2{k_1}} \right)\bar \zeta \left( {2{k_2}} \right)} .
\end{align}
\end{cor}
\begin{cor}\label{cor2.9} If $k>1$ is a positive integer number and $\alpha,\beta$ are real numbers such as $\alpha \beta  = {\pi ^2}$, then
\begin{align}\label{2.17}
&\alpha {\beta ^k}\sum\limits_{n = 1}^\infty  {\frac{{\tanh \left( {\frac{{2n - 1}}
{2}\alpha } \right)}}
{{{{\left( {2n - 1} \right)}^{2k - 1}}}}}  + {\left( { - 1} \right)^k}{\alpha ^k}\beta \sum\limits_{n = 1}^\infty  {\frac{{\tanh \left( {\frac{{2n - 1}}
{2}\beta } \right)}}
{{{{\left( {2n - 1} \right)}^{2k - 1}}}}}\nonumber \\
& =\sum\limits_{\scriptstyle {k_1} + {k_2} = k,\hfill \atop \scriptstyle {k_1}, {k_2} \geqslant 1\hfill}{\frac{{\left( {{2^{2{k_1}}} - 1} \right)\left( {{2^{2{k_2}}} - 1} \right)\left| {{B_{2{k_1}}}} \right|{B_{2{k_2}}}}}
{{\left( {2{k_1}} \right)!\left( {2{k_2}} \right)!}}{\alpha ^{{k_2}}}{\beta ^{{k_1}}}{\pi ^{2k}}} .
\end{align}
\end{cor}
\begin{cor}\label{cor2.10} If $k$ is a non-negative integer number and $\alpha,\beta$ are real numbers such as $\alpha \beta  = {\pi ^2}$, then
\begin{align}\label{2.18}
&{\beta ^k}\sum\limits_{n = 1}^\infty  {\frac{{{{\left( { - 1} \right)}^{n - 1}}}}
{{{{\left( {2n - 1} \right)}^{2k + 1}}\cosh \left( {\frac{{2n - 1}}
{2}\alpha } \right)}}}  + {\left( { - 1} \right)^k}{\alpha ^k}\sum\limits_{n = 1}^\infty  {\frac{{{{\left( { - 1} \right)}^{n - 1}}}}
{{{{\left( {2n - 1} \right)}^{2k + 1}}\cosh \left( {\frac{{2n - 1}}
{2}\beta } \right)}}}\nonumber \\
& =\sum\limits_{\scriptstyle {k_1} + {k_2} = k,\hfill \atop \scriptstyle {k_1}, {k_2} \geqslant 0\hfill}{\frac{{\left| {{E_{2{k_1}}}} \right|{E_{2{k_2}}}}}
{{\left( {2{k_1}} \right)!\left( {2{k_2}} \right)!}}{\alpha ^{{k_2}}}{\beta ^{{k_1}}}{{\left( {\frac{\pi }
{2}} \right)}^{2k + 2}}}.
\end{align}
\end{cor}
Note that the Corollary \ref{cor2.7}, \ref{cor2.9} and \ref{cor2.10} have been given by Chapter 14 of Berndt's \cite{B1989}. It should be emphasized that the Chapter 14 also contains many other types of results.

Then from the definitions of hyperbolic cotangent and tangent, we know that
\[\coth x = 1 + \frac{2}
{{{e^{2x}} - 1}}\quad{\rm and}\quad\tanh x = 1 - \frac{2}
{{{e^{2x}} + 1}}.\]
Hence, substituting above formulae into (\ref{2.14}) and (\ref{2.17}), respectively, by a direct calculation, we obtain the Ramanujan's formula for $\zeta(2n+1)$ and a similar formula that
\begin{align}\label{2.19}
&{\left( {4\beta } \right)^{ - \left( {k - 1} \right)}}\left\{ {\frac{1}
{2}\zeta \left( {2k - 1} \right) + \sum\limits_{n = 1}^\infty  {\frac{1}
{{{n^{2k - 1}}\left( {{e^{2n\alpha }} - 1} \right)}}} } \right\}\nonumber\\
& \quad- {\left( { - 4\alpha } \right)^{ - \left( {k - 1} \right)}}\left\{ {\frac{1}
{2}\zeta \left( {2k - 1} \right) + \sum\limits_{n = 1}^\infty  {\frac{1}
{{{n^{2k - 1}}\left( {{e^{2n\beta }} - 1} \right)}}} } \right\}\nonumber\\
& = \sum\limits_{j = 0}^k {{{\left( { - 1} \right)}^{j - 1}}\frac{{{B_{2j}}{B_{2k - 2j}}}}
{{\left( {2j} \right){\text{!}}\left( {2k - 2j} \right)!}}{\alpha ^j}{\beta ^{k - j}}},
\end{align}
and
\begin{align}\label{2.20}
&{\alpha ^{ - \left( {k - 1} \right)}}\left\{ {\eta \left( {2k - 1} \right) - 2\sum\limits_{n = 1}^\infty  {\frac{1}
{{{e^{\left( {2n - 1} \right)\alpha }} + 1}}} } \right\}\nonumber\\
& \quad- {\left( { - \beta } \right)^{ - \left( {k - 1} \right)}}\left\{ {\eta \left( {2k - 1} \right) - 2\sum\limits_{n = 1}^\infty  {\frac{1}
{{{e^{\left( {2n - 1} \right)\beta }} + 1}}} } \right\}\nonumber\\
&=\sum\limits_{\scriptstyle {k_1} + {k_2} = k,\hfill \atop \scriptstyle {k_1}, {k_2} \geqslant 1\hfill}{\frac{{\left( {{2^{2{k_1}}} - 1} \right)\left( {{2^{2{k_2}}} - 1} \right)\left| {{B_{2{k_1}}}} \right|{B_{2{k_2}}}}}
{{\left( {2{k_1}} \right)!\left( {2{k_2}} \right)!}}{\alpha ^{{k_2}}}{\beta ^{{k_1}}}},
\end{align}
where the function $\eta (s)$ defined by
\[\eta \left( s \right): = \sum\limits_{n = 1}^\infty  {\frac{1}
{{{{\left( {2n - 1} \right)}^s}}}}, \Re(s)>1. \]
It is obvious that
\[\eta \left( s \right) = \left( {1 - {2^{ - s}}} \right)\zeta \left( s \right).\]
The Ramanujan's formula for $\zeta(2n+1)$ (\ref{2.19}) appears as Entry 21(i) in Chapter 14 of Ramanujan's second notebook \cite{R2012}. It also appears in a formerly unpublished manuscript of Ramanujan that was published in its original handwritten form with his lost notebook \cite{R1988}. The first published proof of (\ref{2.19}) is due to S.L. Malurkar \cite{M1925} in the Journal of the Indian Mathematical Society in 1925-1926. This partial manuscript was initially examined in detail by the B.C. Berndt in \cite{B2004}, and by G.E. Andrews and B.C. Berndt in their fourth book on Ramanujan's lost notebook \cite{AB2013}. There exist a huge number of generalizations and analogues of Ramanujan's formula
for $\zeta(2n+1)$. Some of these are discussed in Berndt's paper \cite{B1978} and his book \cite{B1989}. In fact, The Ramanujan's formula for $\zeta(2n+1)$ is a special instance of a general transformation formula for generalized analytic Eisenstein series or, in another formulation, for a vast generalization of the logarithm of the Dedekind eta function, see \cite{BS2017}.

\subsection{Proofs of Theorem \ref{thm2.1}-\ref{thm2.5}}
In the context of this paper, the results of Theorem \ref{thm2.1}-\ref{thm2.5} can be proved by applying the following kernels and base functions (see Table 2.1)
\begin{table}[htbp]\centering
 \begin{tabular}{lll}
  \hline
  kernel function $\xi \left( z\right)$ &  base function $r\left(z\right)$ & combined function $r\left(z\right)\xi \left( z\right)$ \\
 \hline
$\pi \cot \left( {\pi xz} \right)$ & $\frac{{\pi \coth \left( {\pi yz} \right)}}
{{{z^p}}}$ & $\frac{{\pi \cot \left( {\pi xz} \right)\pi \coth \left( {\pi yz} \right)}}
{{{z^p}}}$ \\
$\frac{\pi }
{{\sin \left( {\pi xz} \right)}}$ & $\frac{\pi }
{{{z^p}\sinh \left( {\pi yz} \right)}}$ &$\frac{{{\pi ^2}}}
{{{z^p}\sin \left( {\pi xz} \right)\sinh \left( {\pi yz} \right)}}$\\
$\pi \cot \left( {\pi xz} \right)$ & $\frac{\pi }
{{{z^p}\sinh \left( {\pi yz} \right)}}$& $\frac{{{\pi ^2}\cot \left( {\pi xz} \right)}}
{{{z^p}\sinh \left( {\pi yz} \right)}}$\\
$\pi \tan \left( {\pi xz} \right)$&$\frac{{\pi \tanh \left( {\pi yz} \right)}}
{{{z^p}}}$&$\frac{{{\pi ^2}\tan \left( {\pi xz} \right)\tanh \left( {\pi yz} \right)}}
{{{z^p}}}$\\
$\frac{\pi }
{{\cos \left( {\pi xz} \right)}}$&$\frac{\pi }
{{{z^p}\cosh \left( {\pi yz} \right)}}$&$\frac{{{\pi ^2}}}
{{{z^p}\cos \left( {\pi xz} \right)\cosh \left( {\pi yz} \right)}}$\\
  \hline
 \end{tabular}
 \begin{center}
  \textbf{\footnotesize\bf TABLE 2.1.}\ \ kernels and base functions
  \end{center}
\end{table}\\
To prove Theorem \ref{thm2.1}, we consider the function
\[{f_1}\left( z \right) = \frac{{\pi \cot \left( {\pi xz} \right)\pi \coth \left( {\pi yz} \right)}}
{{{z^p}}},\quad p\geq 2.\]
The only singularities are poles at the $z=\frac {n}{x}$ and $\frac {n}{y}i$, where $n\in \Z=\{0,\pm1,\pm2,\ldots\}$. For $n\in \Z\setminus \{0\}=\{\pm1,\pm2,\pm3,\ldots\}$, these poles are simple, then using the formulae (\ref{2.1}) and (\ref{2.3}), we can deduce that these residues are
\begin{align*}
&{\text{Res}}\left[ {{f_1}\left( z \right),z = \frac{n}
{x}} \right] = {x^{p - 1}}\frac{{\pi \coth \left( {\frac{{\pi ny}}
{x}} \right)}}
{{{n^p}}}, \quad n\in\N\hfill \\
&{\text{Res}}\left[ {{f_1}\left( z \right),z =  - \frac{n}
{x}} \right] = {x^{p - 1}}\frac{{\pi \coth \left( {\frac{{\pi ny}}
{x}} \right)}}
{{{n^p}}}{\left( { - 1} \right)^{p + 1}}, \quad n\in\N\hfill \\
&{\text{Res}}\left[ {{f_1}\left( z \right),z = \frac{n}
{y}i} \right] = {y^{p - 1}}\frac{{\pi \coth \left( {\frac{{\pi nx}}
{y}} \right)}}
{{{n^p}}}{\left( { - 1} \right)^{p + 1}}{i^{p + 1}}, \quad n\in\N\hfill \\
&{\text{Res}}\left[ {{f_1}\left( z \right),z =  - \frac{n}
{y}i} \right] = {y^{p - 1}}\frac{{\pi \coth \left( {\frac{{\pi nx}}
{y}} \right)}}
{{{n^p}}}{i^{p + 1}}, \quad n\in\N.
\end{align*}
When $z\Rightarrow 0$, using the formulae (\ref{2.1}) and (\ref{2.3}) again, we can find that
\begin{align}\label{2.21}{f_1}\left( z \right)\mathop  = \limits^{z \to 0} \frac{1}
{{{z^p}xy}}\left\{ \begin{gathered}
  \frac{1}
{{{z^2}}} - 2\sum\limits_{k = 1}^\infty  {{{\left( { - 1} \right)}^k}{y^{2k}}\zeta \left( {2k} \right){z^{2k - 2}}}  - 2\sum\limits_{k = 1}^\infty  {{x^{2k}}\zeta \left( {2k} \right){z^{2k - 2}}}  \hfill \\
  \;\;\;\; + 4\sum\limits_{{k_1},{k_2} = 1}^\infty  {{{\left( { - 1} \right)}^{{k_2}}}{x^{2{k_1}}}{y^{2{k_2}}}\zeta \left( {2{k_1}} \right)\zeta \left( {2{k_2}} \right){z^{2{k_1} + 2{k_2} - 2}}}  \hfill \\
\end{gathered}  \right\}.
\end{align}
Hence, the residue of the pole of order $p+2$ at $0$ is found to be
\begin{align}\label{2.22}
{\text{Res}}\left[ {{f_1}\left( z \right),z = 0} \right] &= \frac{1}
{{\left( {p + 1} \right)!}}\frac{{{d^{p + 1}}}}
{{d{z^{p + 1}}}}{\left. {\left( {{z^{p + 2}}{f_1}\left( z \right)} \right)} \right|_{z = 0}}\nonumber\\
& = \frac{1}
{{xy}}\left\{ \begin{gathered}
   - 2\sum\limits_{2k = p + 1,k \geqslant 1} {{{\left( { - 1} \right)}^k}{y^{2k}}\zeta \left( {2k} \right)}  - 2\sum\limits_{2k = p + 1,k \geqslant 1} {{x^{2k}}\zeta \left( {2k} \right)}  \hfill \\
   + 4\sum\limits_{2{k_1} + 2{k_2} = p + 1,{k_1},{k_2} \geqslant 1} {{{\left( { - 1} \right)}^{{k_2}}}{x^{2{k_1}}}{y^{2{k_2}}}\zeta \left( {2{k_1}} \right)\zeta \left( {2{k_2}} \right)}  \hfill \\
\end{gathered}  \right\}.
\end{align}
Applying Lemma \ref{lem1} and summing these residues with setting $p=2k-1\ (k>1)$ yields the statement of the Theorem \ref{thm2.1}.\hfill$\square$

Similarly, to prove the Theorem \ref{thm2.2}-\ref{thm2.5}, we need to consider the following four combined functions
\begin{align*}
&{f_2}\left( z \right) = \frac{{{\pi ^2}}}
{{{z^p}\sin \left( {\pi xz} \right)\sinh \left( {\pi yz} \right)}},\quad{f_3}\left( z \right) = \frac{{{\pi ^2}\cot \left( {\pi xz} \right)}}
{{{z^p}\sinh \left( {\pi yz} \right)}},\\
&{f_4}\left( z \right) = \frac{{{\pi ^2}\tan \left( {\pi xz} \right)\tanh \left( {\pi yz} \right)}}
{{{z^p}}},\quad{f_5}\left( z \right) = \frac{{{\pi ^2}}}
{{{z^p}\cos \left( {\pi xz} \right)\cosh \left( {\pi yz} \right)}}.
\end{align*}
Then applying the same arguments as in the proof of Theorem \ref{thm2.1} with the help of formulae (\ref{2.1})-(\ref{2.8}) we may easily deduce these results. Note that in prove Theorem \ref{thm2.4} and \ref{thm2.5}, we used the following results of Power series expansions \cite{GR2007}
\begin{align*}
&\tan x = \sum\limits_{n = 1}^\infty  {\frac{{{2^{2n}}\left( {{2^{2n}} - 1} \right)}}
{{\left( {2n} \right)!}}\left| {{B_{2n}}} \right|{x^{2n - 1}}} ,{x^2} < \frac{{{\pi ^2}}}
{4}, \hfill \\
&\tanh  = \sum\limits_{n = 1}^\infty  {\frac{{{2^{2n}}\left( {{2^{2n}} - 1} \right)}}
{{\left( {2n} \right)!}}{B_{2n}}{x^{2n - 1}}} ,{x^2} < \frac{{{\pi ^2}}}
{4}, \hfill \\
&\sec x = \sum\limits_{n = 0}^\infty  {\frac{{\left| {{E_{2n}}} \right|}}
{{\left( {2n} \right)!}}{x^{2n}}} ,{x^2} < \frac{{{\pi ^2}}}
{4}, \hfill \\
&\operatorname{sech} x = \sum\limits_{n = 0}^\infty  {\frac{{{E_{2n}}}}
{{\left( {2n} \right)!}}{x^{2n}}} ,{x^2} < \frac{{{\pi ^2}}}
{4}.
\end{align*}
\subsection{Some simple examples}
From subsection 2.1, we can obtain many results involving hyperbolic functions. Some interesting (known or
new) examples follows:
\begin{align}
& \sum\limits_{n = 1}^\infty  {\frac{{\pi \coth \left( {\pi n} \right)}}
{{{n^3}}}}  = \frac{7}
{{180}}{\pi ^4}, \hfill \\
 & \sum\limits_{n = 1}^\infty  {\frac{{\pi {{\left( { - 1} \right)}^{n - 1}}}}
{{{n^3}\sinh \left( {\pi n} \right)}}}  = \frac{1}
{{360}}{\pi ^4}, \hfill \\
 & \sum\limits_{n = 1}^\infty  {\frac{{\pi \tanh \left( {\frac{{2n - 1}}
{2}\pi } \right)}}
{{{{\left( {2n - 1} \right)}^3}}}}  = \frac{{{\pi ^4}}}
{{32}}, \hfill \\
 & \sum\limits_{n = 1}^\infty  {\frac{{{{\left( { - 1} \right)}^{n - 1}}}}
{{\left( {2n - 1} \right)\cosh \left( {\frac{{2n - 1}}
{2}\pi } \right)}}}  = \frac{\pi }
{8},  \\
&\sum\limits_{n = 1}^\infty  {\frac{{\pi \coth \left( {\pi n} \right)}}
{{{n^{4k - 1}}}}}  = \zeta \left( {4k} \right) - \sum\limits_{\scriptstyle {k_1} + {k_2} =2 k,\hfill \atop \scriptstyle {k_1}, {k_2} \geqslant 1\hfill} {{{\left( { - 1} \right)}^{{k_2}}}\zeta \left( {2{k_1}} \right)\zeta \left( {2{k_2}} \right)} ,\\
&\sum\limits_{n = 1}^\infty  {\frac{{\pi {{\left( { - 1} \right)}^{n - 1}}}}
{{{n^{4k - 1}}\sinh \left( {\pi n} \right)}}}  = \bar \zeta \left( {4k} \right) + \sum\limits_{\scriptstyle {k_1} + {k_2} =2 k,\hfill \atop \scriptstyle {k_1}, {k_2} \geqslant 1\hfill} {{{\left( { - 1} \right)}^{{k_2}}}\bar \zeta \left( {2{k_1}} \right)\bar \zeta \left( {2{k_2}} \right)} ,\\
&\sum\limits_{n = 1}^\infty  {\frac{{\pi \tanh \left( {\frac{{2n - 1}}
{2}\pi } \right)}}
{{{{\left( {2n - 1} \right)}^{4k - 1}}}}}  =\sum\limits_{\scriptstyle {k_1} + {k_2} =2 k,\hfill \atop \scriptstyle {k_1}, {k_2} \geqslant 1\hfill} {\frac{{\left( {{2^{2{k_1}}} - 1} \right)\left( {{2^{2{k_2}}} - 1} \right)\left| {{B_{2{k_1}}}} \right|{B_{2{k_2}}}}}
{{2\left( {2{k_1}} \right)!\left( {2{k_2}} \right)!}}{\pi ^{4k}}} ,\\
&\sum\limits_{n = 1}^\infty  {\frac{{\pi {{\left( { - 1} \right)}^{n - 1}}}}
{{{{\left( {2n - 1} \right)}^{4k + 1}}\cosh \left( {\frac{{2n - 1}}
{2}\pi } \right)}}}  = \sum\limits_{\scriptstyle {k_1} + {k_2} =2 k,\hfill \atop \scriptstyle {k_1}, {k_2} \geqslant 0\hfill}{\frac{{\left| {{E_{2{k_1}}}} \right|{E_{2{k_2}}}}}
{{2\left( {2{k_1}} \right)!\left( {2{k_2}} \right)!}}{{\left( {\frac{\pi }
{2}} \right)}^{4k + 2}}},\\
&\sum\limits_{n = 1}^\infty  {\frac{\pi }
{{{n^{2k - 1}}\sinh \left( {\pi n} \right)}}}  + {\left( { - 1} \right)^k}\sum\limits_{n = 1}^\infty  {\frac{{\pi \coth \left( {\pi n} \right)}}
{{{n^{2k - 1}}}}{{\left( { - 1} \right)}^n}} \nonumber\\
 &\quad\quad\quad\quad\quad\quad\quad\quad\quad\quad\quad= \zeta \left( {2k} \right) - {\left( { - 1} \right)^k}\bar \zeta \left( {2k} \right) + 2\sum\limits_{{k_1} + {k_2} = k} {{{\left( { - 1} \right)}^{{k_2}}}\zeta \left( {2{k_1}} \right)\bar \zeta \left( {2{k_2}} \right)}  ,\\
 &\sum\limits_{n = 1}^\infty  {\frac{\pi }
{{n\sinh \left( {\frac{{\pi ny}}
{x}} \right)}}}  - \sum\limits_{n = 1}^\infty  {\frac{{\pi \coth \left( {\frac{{\pi nx}}
{y}} \right)}}
{n}{{\left( { - 1} \right)}^n}}  = \frac{1}
{2}\zeta \left( 2 \right)\left( {2x + y} \right),\\
 &\sum\limits_{n = 1}^\infty  {\left\{ {\frac{1}
{{\left( {2n - 1} \right)\cosh \left( {\frac{{2n - 1}}
{{2x}}\pi y} \right)}} + \frac{1}
{{\left( {2n - 1} \right)\cosh \left( {\frac{{2n - 1}}
{{2y}}\pi x} \right)}}} \right\}{{\left( { - 1} \right)}^{n - 1}}}  = \frac{\pi }
{4},\\
&{x^2}\sum\limits_{n = 1}^\infty  {\frac{{\pi \tanh \left( {\frac{{2n - 1}}
{{2x}}\pi y} \right)}}
{{{{\left( {2n - 1} \right)}^3}}}}  + {y^2}\sum\limits_{n = 1}^\infty  {\frac{{\pi \tanh \left( {\frac{{2n - 1}}
{{2y}}\pi x} \right)}}
{{{{\left( {2n - 1} \right)}^3}}}}  = \frac{{{\pi ^4}}}
{{16}}xy,\\
&{x^2}\sum\limits_{n = 1}^\infty  {\frac{{\pi {{\left( { - 1} \right)}^{n - 1}}}}
{{{{\left( {2n - 1} \right)}^3}\cosh \left( {\frac{{2n - 1}}
{{2x}}\pi y} \right)}}}  - {y^2}\sum\limits_{n = 1}^\infty  {\frac{{\pi {{\left( { - 1} \right)}^{n - 1}}}}
{{{{\left( {2n - 1} \right)}^3}\cosh \left( {\frac{{2n - 1}}
{{2y}}\pi x} \right)}}}  = \frac{{{\pi ^4}}}
{{32}}\left( {{x^2} - {y^2}} \right),
\end{align}
\section{The partial fraction expansion of several quadratic products involving hyperbolic and trigonometric function}
In his paper \cite{R1916}, Ramanujan's mistake arose when he attempted to find the partial fraction expansion of $\cot \left( {\sqrt {w\alpha } } \right)\coth \left( {\sqrt {w\beta } } \right)$ that
\begin{align}\label{3.1}
\frac{\pi }
{2}\cot \left( {\sqrt {w\alpha } } \right)\coth \left( {\sqrt {w\beta } } \right) = \frac{1}
{{2w}} + \sum\limits_{n = 1}^\infty  {\left\{ {\frac{{n\alpha \coth \left( {n\alpha } \right)}}
{{w + {n^2}\alpha }} + \frac{{n\beta \coth \left( {n\beta } \right)}}
{{w - {n^2}\beta }}} \right\}}\quad (\alpha \beta  = {\pi ^2}).
\end{align}
B.C. Berndt and A. Straub \cite{BS2017} apply the Mittag-Leffler Theorem to find the mistake. Moreover, they offer a corrected version of (\ref{3.1}) that
\begin{align}\label{3.2}
\frac{\pi }
{2}\cot \left( {\sqrt {w\alpha } } \right)\coth \left( {\sqrt {w\beta } } \right) = \frac{1}
{{2w}} + \frac{1}
{2}\log \frac{\beta }
{\alpha } + \sum\limits_{n = 1}^\infty  {\left\{ {\frac{{n\alpha \coth \left( {n\alpha } \right)}}
{{w + {n^2}\alpha }} + \frac{{n\beta \coth \left( {n\beta } \right)}}
{{w - {n^2}\beta }}} \right\}}
\end{align}
by using the correct partial fraction decomposition given by R. Sitaramachandrarao \cite{S}
\begin{align}\label{3.3}
{\pi ^2}xy\cot \left( {\pi x} \right)\coth \left( {\pi y} \right) = 1 + 2\zeta \left( 2 \right)\left( {{y^2} - {x^2}} \right) - 2\pi xy\sum\limits_{n = 1}^\infty  {\left\{ {\frac{{{x^2}\coth \left( {\frac{{\pi ny}}
{x}} \right)}}
{{n\left( {{n^2} - {x^2}} \right)}} + \frac{{{y^2}\coth \left( {\frac{{\pi nx}}
{y}} \right)}}
{{n\left( {{n^2} - {y^2}} \right)}}} \right\}}.
\end{align}
In fact, if we consider the function
\begin{align}\label{3.4}
f\left( z \right) = \frac{{\pi^2 \cot \left( {\pi xz} \right)\coth \left( {\pi yz} \right)}}
{{z - a}}\quad {\rm or}\quad \frac{{\pi^2 \cot \left( {\pi xz} \right)\coth \left( {\pi yz} \right)z}}
{{{z^2} - {a^2}}}.
\end{align}
Then applying the same arguments as in the proof of Theorem \ref{thm2.1}, we may deduce the formula (\ref{3.1}). However, from Lemma \ref{lem1} and the definition of hyperbolic cotangent, we find that
\begin{align*}
\frac{{\coth \left( {\pi yz} \right)}}
{{z - a}} &= \frac{1}
{{z - a}}\left( {1 + \frac{2}
{{{e^{2\pi yz}} - 1}}} \right)\\
& = \frac{1}
{{z - a}} + \frac{2}
{{z - a\left( {{e^{2\pi yz}} - 1} \right)}}.
\end{align*}
Hence, if $z\rightarrow \infty$, we have
\[\frac{{\coth \left( {\pi yz} \right)}}
{{z - a}} = O\left( {{z^{ - 1}}} \right) \ne O\left( {{z^{ - 2}}} \right).\]
Therefore, according to the condition of Lemma \ref{lem1} holds, we know that we can't use the Lemma \ref{lem1} to evaluate the residues of function $f(z)$ defined by (\ref{3.4}). But if we consider the function
\[g_1\left( z \right) = \frac{{\pi^2 \cot \left( {\pi xz} \right)\coth \left( {\pi yz} \right)}}
{{z\left( {{z^2} - {a^2}} \right)}}.\]
Then by a similar argument as in the proof of Theorem \ref{thm2.1}, we deduce that
\begin{align*}
 & {\text{Res}}\left[ {g_1\left( z \right),z =  \pm \frac{n}
{x}} \right] = {x^2}\frac{{\pi \coth \left( {\frac{{\pi ny}}
{x}} \right)}}
{{n\left( {{n^2} - {a^2}{x^2}} \right)}}, \hfill \\
&  {\text{Res}}\left[ {g_1\left( z \right),z =  \pm \frac{n}
{y}i} \right] = {y^2}\frac{{\pi \coth \left( {\frac{{\pi nx}}
{y}} \right)}}
{{n\left( {{n^2} + {a^2}{y^2}} \right)}}, \hfill \\
 & {\text{Res}}\left[ {g_1\left( z \right),z =  \pm a} \right] = \frac{{{\pi ^2}\cot \left( {\pi xa} \right)\coth \left( {\pi ya} \right)}}
{{2{a^2}}}, \hfill \\
 & {\text{Res}}\left[ {g_1\left( z \right),z = 0} \right] =  - \frac{1}
{{{a^4}xy}} - \frac{{2\zeta \left( 2 \right)}}
{{{a^2}}}\left( {\frac{y}
{x} - \frac{x}
{y}} \right).
\end{align*}
Hence, summing all residues yields the formula
\begin{align}\label{3.5}
{\pi ^2}xy\cot \left( {\pi xa} \right)\coth \left( {\pi ya} \right) =& \frac{1}
{{{a^2}}} + 2\zeta \left( 2 \right)\left( {{y^2} - {x^2}} \right)\nonumber \\&- 2\pi xy\sum\limits_{n = 1}^\infty  {\left\{ {\frac{{{a^2}{x^2}\coth \left( {\frac{{\pi ny}}
{x}} \right)}}
{{n\left( {{n^2} - {a^2}{x^2}} \right)}} + \frac{{{a^2}{y^2}\coth \left( {\frac{{\pi nx}}
{y}} \right)}}
{{n\left( {{n^2} - {a^2}{y^2}} \right)}}} \right\}}.
\end{align}
Letting $a=1$ in (\ref{3.5}), the result is (\ref{3.3}).

In this section, by using the residue theorem, we give the partial fraction expansion of several quadratic products involving hyperbolic and trigonometric function of the form
\[\frac{\pi }
{{\sin \left( {\sqrt {w\alpha } } \right)\sinh \left( {\sqrt {w\beta } } \right)}},\frac{{\pi \cot \left( {\sqrt {w\alpha } } \right)}}
{{\sinh \left( {\sqrt {w\beta } } \right)}},\pi \tan \left( {\sqrt {w\alpha } } \right)\tanh \left( {\sqrt {w\beta } } \right),\frac{\pi }
{{\cos \left( {\sqrt {w\alpha } } \right)\cosh \left( {\sqrt {w\beta } } \right)}},\]
where $\alpha \beta  = {\pi ^2}$.
\subsection{Some identities on series of hyperbolic function}
Applying Lemma \ref{lem1} and formulae (\ref{2.1})-(\ref{2.8}), we can get the following theorem.
\begin{thm} For any reals $a,x,y\neq 0$ with $xa \ne \pm\frac{{2n - 1}}{2}\ (n\in\Z)$, then the following identities hold:
\begin{align}\label{3.6}
\frac{{{\pi ^2}}}
{{\sin \left( {\pi xa} \right)\sinh \left( {\pi ya} \right)}} = \frac{1}
{{{a^2}xy}} + 2\pi \sum\limits_{n = 1}^\infty  {\left\{ {\frac{{n{{\left( { - 1} \right)}^n}}}
{{\left( {{n^2} + {a^2}{y^2}} \right)\sinh \left( {\frac{{\pi nx}}
{y}} \right)}} - \frac{{n{{\left( { - 1} \right)}^n}}}
{{\left( {{n^2} - {a^2}{x^2}} \right)\sinh \left( {\frac{{\pi ny}}
{x}} \right)}}} \right\}},
\end{align}

\begin{align}\label{3.7}
\frac{{{\pi ^2}\cot \left( {\pi xa} \right)}}
{{\sinh \left( {\pi ya} \right)}} = \frac{1}
{{{a^2}xy}} + 2\pi \sum\limits_{n = 1}^\infty  {\left\{ {\frac{{\pi n\coth \left( {\frac{{\pi nx}}
{y}} \right)}}
{{\left( {{n^2} + {a^2}{y^2}} \right)}} (-1)^n- \frac{n}
{{\left( {{n^2} - {a^2}{x^2}} \right)\sinh \left( {\frac{{\pi ny}}
{x}} \right)}}} \right\}} ,
\end{align}

\begin{align}\label{3.8}
\frac{{{\pi ^2}\tan \left( {\pi xa} \right)\tanh \left( {\pi ya} \right)}}
{{16{a^2}}} = \sum\limits_{n = 1}^\infty  {\left\{ {\frac{{\pi {x^2}\tanh \left( {\frac{{2n - 1}}
{{2x}}\pi y} \right)}}
{{\left( {2n - 1} \right)\left[ {{{\left( {2n - 1} \right)}^2} - {{\left( {2xa} \right)}^2}} \right]}} + \frac{{\pi {y^2}\tanh \left( {\frac{{2n - 1}}
{{2y}}\pi x} \right)}}
{{\left( {2n - 1} \right)\left[ {{{\left( {2n - 1} \right)}^2} + {{\left( {2ya} \right)}^2}} \right]}}} \right\}} ,
\end{align}

\begin{align}\label{3.9}
\frac{{{\pi ^2}}}
{{16{a^2}\cos \left( {\pi xa} \right)\cosh \left( {\pi ya} \right)}} =& \frac{{{\pi ^2}}}
{{16{a^2}}} + \sum\limits_{n = 1}^\infty  {\frac{{\pi {y^2}{{\left( { - 1} \right)}^n}}}
{{\left( {2n - 1} \right)\left[ {{{\left( {2n - 1} \right)}^2} + {{\left( {2ya} \right)}^2}} \right]\cosh \left( {\frac{{2n - 1}}
{{2y}}\pi x} \right)}}} \nonumber\\
& - \sum\limits_{n = 1}^\infty  {\frac{{\pi {x^2}{{\left( { - 1} \right)}^n}}}
{{\left( {2n - 1} \right)\left[ {{{\left( {2n - 1} \right)}^2} - {{\left( {2xa} \right)}^2}} \right]\cosh \left( {\frac{{2n - 1}}
{{2x}}\pi y} \right)}}} .
\end{align}
\end{thm}
\pf To prove the formulae (\ref{3.6})-(\ref{3.9}), we consider the four functions (see Table \ref{t2})
\begin{table}[htbp]\centering\label{t2}
 \begin{tabular}{lll}
  \hline
  kernel function $\xi \left( z\right)$ &  base function $r\left(z\right)$ & combined function $r\left(z\right)\xi \left( z\right)$ \\
$\frac{\pi }
{{\sin \left( {\pi xz} \right)}}$ & $\frac{\pi }
{{{(z-a)}\sinh \left( {\pi yz} \right)}}$ &$\frac{{{\pi ^2}}}
{{{(z-a)}\sin \left( {\pi xz} \right)\sinh \left( {\pi yz} \right)}}$\\
$\pi \cot \left( {\pi xz} \right)$ & $\frac{\pi }
{{{(z-a)}\sinh \left( {\pi yz} \right)}}$& $\frac{{{\pi ^2}\cot \left( {\pi xz} \right)}}
{{{(z-a)}\sinh \left( {\pi yz} \right)}}$\\
$\pi \tan \left( {\pi xz} \right)$&$\frac{{\pi \tanh \left( {\pi yz} \right)}}
{{{z(z^2-a^2)}}}$&$\frac{{{\pi ^2}\tan \left( {\pi xz} \right)\tanh \left( {\pi yz} \right)}}
{{{z(z^2-a^2)}}}$\\
$\frac{\pi }
{{\cos \left( {\pi xz} \right)}}$&$\frac{\pi }
{{{z(z^2-a^2)}\cosh \left( {\pi yz} \right)}}$&$\frac{{{\pi ^2}}}
{{{z(z^2-a^2)}\cos \left( {\pi xz} \right)\cosh \left( {\pi yz} \right)}}$\\
  \hline
 \end{tabular}
 \begin{center}
  \textbf{\footnotesize\bf TABLE 3.1.}\ \ kernels functions and base functions of parametric $a$
  \end{center}
\end{table}
\begin{align*}
& {g_2}\left( z \right) = \frac{{{\pi ^2}}}
{{\left( {z - a} \right)\sin \left( {\pi xz} \right)\sinh \left( {\pi yz} \right)}}, \hfill \\
&  {g_3}\left( z \right) = \frac{{{\pi ^2}\cot \left( {\pi xz} \right)}}
{{\left( {z - a} \right)\sinh \left( {\pi yz} \right)}}, \hfill \\
 & {g_4}\left( z \right) = \frac{{{\pi ^2}\tan \left( {\pi xz} \right)\tanh \left( {\pi yz} \right)}}
{{z\left( {{z^2} - {a^2}} \right)}}, \hfill \\
 & {g_5}\left( z \right) = \frac{{{\pi ^2}}}
{{z\left( {{z^2} - {a^2}} \right)\cos \left( {\pi xz} \right)\cosh \left( {\pi yz} \right)}}.
\end{align*}
Then by a similar argument as in the proof of formula (\ref{3.5}), we may easily deduce these results (\ref{3.6})-(\ref{3.9}).\hfill$\square$

If $a=1$ these results were also proved in Berndt's book \cite{B1989} by another method. A large number of formulae aries from partial fraction decompositions.  For
more details, see Chapter 14, pp 240-299 of his book \cite{B1989}. It is obvious that proceeding in a similar method to evaluation of the (\ref{3.6})-(\ref{3.9}), it is possible
that some of similar results of partial fraction expansion can be obtained. For example, we can also deduce the following well known identities (\cite{B1989})
\begin{align*}
&\sum\limits_{n = 1}^\infty  {\frac{{n{{\left( { - 1} \right)}^n}}}
{{\left( {{n^4} - {a^4}} \right)\sinh \left( {\pi n} \right)}}}  = \frac{1}
{{4{a^2}\pi }}\left( {\frac{1}
{{{a^2}}} - \frac{{{\pi ^2}}}
{{\sin \left( {\pi a} \right)\sinh \left( {\pi a} \right)}}} \right), \hfill \\
&\sum\limits_{n = 1}^\infty  {\frac{{n\coth \left( {\pi n} \right)}}
{{{n^4} - {a^4}}}}  = \frac{1}
{{4{a^2}\pi }}\left( {\frac{1}
{{{a^2}}} - {\pi ^2}\cot \left( {\pi a} \right)\coth \left( {\pi a} \right)} \right).
\end{align*}
\subsection{Four partial fraction expansion formulae}
Setting the transformation $\alpha  = \pi \frac{x}{y},\beta  = \pi \frac{y}{x},w = \pi xy{a^2}$ in formulae (\ref{3.6})-(\ref{3.9}), by a simple calculation we can obtain these partial fraction expansions
\begin{align}\label{3.10}
\frac{\pi }
{{2\sin \left( {\sqrt {w\alpha } } \right)\sinh \left( {\sqrt {w\beta } } \right)}} = \frac{1}
{{2w}} + \sum\limits_{n = 1}^\infty  {\left\{ {\frac{{\beta n{{\left( { - 1} \right)}^n}}}
{{\left( {w - \beta {n^2}} \right)\sinh \left( {\beta n} \right)}} + \frac{{n{{\left( { - 1} \right)}^n}}}
{{\left( {w + \alpha {n^2}} \right)\sinh \left( {\alpha n} \right)}}} \right\}},
\end{align}

\begin{align}\label{3.11}
\frac{{\pi \cot \left( {\sqrt {w\alpha } } \right)}}
{{2\sinh \left( {\sqrt {w\beta } } \right)}} = \frac{1}
{{2w}} + \sum\limits_{n = 1}^\infty  {\left\{ {\frac{{\beta n}}
{{\left( {w - \beta {n^2}} \right)\sinh \left( {\beta n} \right)}} + \frac{{\alpha n\coth \left( {\alpha n} \right)}}
{{w + \alpha {n^2}}}{{\left( { - 1} \right)}^n}} \right\}} ,
\end{align}

\begin{align}\label{3.12}
\frac{{\pi \tan \left( {\sqrt {w\alpha } } \right)\tanh \left( {\sqrt {w\beta } } \right)}}
{{16w}} = \sum\limits_{n = 1}^\infty  {\left\{ {\frac{{\tanh \left( {\frac{{2n - 1}}
{2}\beta } \right)}}
{{\left( {2n - 1} \right)\left[ {\beta {{\left( {2n - 1} \right)}^2} - 4w} \right]}} + \frac{{\tanh \left( {\frac{{2n - 1}}
{2}\alpha } \right)}}
{{\left( {2n - 1} \right)\left[ {\alpha {{\left( {2n - 1} \right)}^2} + 4w} \right]}}} \right\}} ,
\end{align}

\begin{align}\label{3.13}
\frac{\pi }
{{16w\cos \left( {\sqrt {w\alpha } } \right)\cosh \left( {\sqrt {w\beta } } \right)}} =& \frac{\pi }
{{16w}} + \sum\limits_{n = 1}^\infty  {\frac{{{{\left( { - 1} \right)}^n}}}
{{\left( {2n - 1} \right)\left[ {4w - \beta {{\left( {2n - 1} \right)}^2}} \right]\cosh \left( {\frac{{2n - 1}}
{2}\beta } \right)}}}\nonumber \\
&+ \sum\limits_{n = 1}^\infty  {\frac{{{{\left( { - 1} \right)}^n}}}
{{\left( {2n - 1} \right)\left[ {4w + \alpha {{\left( {2n - 1} \right)}^2}} \right]\cosh \left( {\frac{{2n - 1}}
{2}\alpha } \right)}}}.
\end{align}
\section{Some results on sums of quadratic hyperbolic functions}
In this section, we will evaluate some series involving quadratic hyperbolic functions, and establish many relations. First, we consider the function
\[{h_1}\left( z \right) = \frac{\pi }
{{\sin \left( {\pi z} \right)}}\frac{{{\pi ^2}{{\coth }^2}\left( {\pi z} \right)}}
{{{z^{2k}}}},\quad k\in\N.\]
The function only singularities are poles at the integers and $\pm ni$. By a direct calculation, we find that
\begin{align*}
&{\text{Res}}\left[ {{h_1}\left( z \right),z =  \pm n} \right] = \frac{{{\pi ^2}{{\coth }^2}\left( {\pi z} \right)}}
{{{n^{2k}}}}{\left( { - 1} \right)^n},\quad n\in\N,\\
&{\text{Res}}\left[ {{h_1}\left( z \right),z =  \pm ni} \right] = 2k{\left( { - 1} \right)^k}\frac{\pi }
{{{n^{2k + 1}}\sinh \left( {\pi n} \right)}} + {\left( { - 1} \right)^k}\frac{{{\pi ^2}\cosh \left( {\pi n} \right)}}
{{{n^{2k}}{{\sinh }^2}\left( {\pi n} \right)}},\quad n\in\N,\\
&{\text{Res}}\left[ {{h_1}\left( z \right),z = 0} \right] = 2\zeta \left( {2k + 2} \right)\left\{ {1 - {2^{ - 2k - 1}} - 2{{\left( { - 1} \right)}^{k + 1}}} \right\}\\
&\quad\quad\quad\quad\quad\quad\quad\quad\quad + 2\sum\limits_{\scriptstyle {k_1} + {k_2} = k,\hfill \atop \scriptstyle {k_1}, {k_2} \geqslant 1\hfill} {\zeta \left( {2{k_1}} \right)\zeta \left( {2{k_2}} \right)\left\{ {2{{\left( { - 1} \right)}^{k + 1}} - 4{{\left( { - 1} \right)}^{{k_2}}}\left( {1 - {2^{1 - 2{k_1}}}} \right)} \right\}} \\
&\quad\quad\quad\quad\quad\quad\quad\quad\quad  + 8\sum\limits_{\scriptstyle {k_1} + {k_2}+k_3 = k,\hfill \atop \scriptstyle {k_1}, {k_2},k_3 \geqslant 1\hfill} {{{\left( { - 1} \right)}^{{k_1} + {k_2}}}\left( {1 - {2^{1 - 2{k_3}}}} \right)\zeta \left( {2{k_1}} \right)\zeta \left( {2{k_2}} \right)\zeta \left( {2{k_3}} \right)}.
\end{align*}
Hence, applying Lemma \ref{lem1} we obtain the result that
\begin{align}\label{4.1}
&\sum\limits_{n = 1}^\infty  {\frac{{{\pi ^2}{{\coth }^2}\left( {\pi z} \right)}}
{{{n^{2k}}}}{{\left( { - 1} \right)}^n}}  + 2k{\left( { - 1} \right)^k}\sum\limits_{n = 1}^\infty  {\frac{\pi }
{{{n^{2k + 1}}\sinh \left( {\pi n} \right)}}}  + {\left( { - 1} \right)^k}\sum\limits_{n = 1}^\infty  {\frac{{{\pi ^2}\cosh \left( {\pi n} \right)}}
{{{n^{2k}}{{\sinh }^2}\left( {\pi n} \right)}}}  \nonumber \\
 &  + \zeta \left( {2k + 2} \right)\left\{ {1 - {2^{ - 2k - 1}} - 2{{\left( { - 1} \right)}^{k + 1}}} \right\} \nonumber \\
  & + \sum\limits_{\scriptstyle {k_1} + {k_2} = k,\hfill \atop \scriptstyle {k_1}, {k_2} \geqslant 1\hfill}  {\zeta \left( {2{k_1}} \right)\zeta \left( {2{k_2}} \right)\left\{ {2{{\left( { - 1} \right)}^{k + 1}} - 4{{\left( { - 1} \right)}^{{k_2}}}\left( {1 - {2^{1 - 2{k_1}}}} \right)} \right\}}  \nonumber\\
  & + 4\sum\limits_{\scriptstyle {k_1} + {k_2}+k_3 = k,\hfill \atop \scriptstyle {k_1}, {k_2},k_3 \geqslant 1\hfill} {{{\left( { - 1} \right)}^{{k_1} + {k_2}}}\left( {1 - {2^{1 - 2{k_3}}}} \right)\zeta \left( {2{k_1}} \right)\zeta \left( {2{k_2}} \right)\zeta \left( {2{k_3}} \right)}  = 0.
\end{align}
As simple example is as follows: setting $k=1$ in (\ref{4.1}), then
\[\sum\limits_{n = 1}^\infty  {\frac{{{\pi ^2}{{\coth }^2}\left( {\pi z} \right)}}
{{{n^2}}}{{\left( { - 1} \right)}^n}}  - 2\sum\limits_{n = 1}^\infty  {\frac{\pi }
{{{n^3}\sinh \left( {\pi n} \right)}}}  - \sum\limits_{n = 1}^\infty  {\frac{{{\pi ^2}\cosh \left( {\pi n} \right)}}
{{{n^2}{{\sinh }^2}\left( {\pi n} \right)}}}  + \frac{{71}}
{8}\zeta \left( 4 \right) = 0.\]
Similarly, by considering the function
\begin{align*}
 & {h_2}\left( z \right) = \frac{\pi }
{{\sin \left( {\pi z} \right)}}\frac{{{\pi ^2}}}
{{{z^{2k}}{{\sinh }^2}\left( {\pi z} \right)}}, \hfill \\
 & {h_3}\left( z \right) = \pi \cot \left( {\pi z} \right)\frac{{{\pi ^2}}}
{{{z^{2k}}{{\sinh }^2}\left( {\pi z} \right)}}, \hfill \\
 & {h_4}\left( z \right) = \pi \cot \left( {\pi z} \right)\frac{{{\pi ^2}{{\coth }^2}\left( {\pi z} \right)}}
{{{z^{2k}}}},
\end{align*}
we can deduce the following results
\begin{align}\label{4.2}
&\sum\limits_{n = 1}^\infty  {\frac{{{\pi ^2}}}
{{{n^{2k}}{{\sinh }^2}\left( {\pi n} \right)}}{{\left( { - 1} \right)}^n}}  + 2k{\left( { - 1} \right)^k}\sum\limits_{n = 1}^\infty  {\frac{\pi }
{{{n^{2k + 1}}\sinh \left( {\pi n} \right)}}}  + {\left( { - 1} \right)^k}\sum\limits_{n = 1}^\infty  {\frac{{{\pi ^2}\cosh \left( {\pi n} \right)}}
{{{n^{2k}}{{\sinh }^2}\left( {\pi n} \right)}}}\nonumber \\
&+ \zeta \left( {2k + 2} \right)\left( {1 - {2^{ - 2k - 1}}} \right)\left( {2 + {{\left( { - 1} \right)}^{k + 1}}} \right) + 2\sum\limits_{\scriptstyle {k_1} + {k_2} = k,\hfill \atop \scriptstyle {k_1}, {k_2} \geqslant 1\hfill} {{{\left( { - 1} \right)}^{{k_2}}}\left( {{{\left( { - 1} \right)}^{{k_1}}} + 2} \right)\bar \zeta \left( {2{k_1}} \right)\bar \zeta \left( {2{k_2}} \right)}\nonumber \\
& + 4\sum\limits_{\scriptstyle {k_1} + {k_2}+k_3 = k,\hfill \atop \scriptstyle {k_1}, {k_2},k_3 \geqslant 1\hfill} {{{\left( { - 1} \right)}^{{k_1} + {k_2}}}\bar \zeta \left( {2{k_1}} \right)\bar \zeta \left( {2{k_2}} \right)\bar \zeta \left( {2{k_3}} \right)}  = 0,
\end{align}
\begin{align}\label{4.3}
 & \left( {1 + {{\left( { - 1} \right)}^k}} \right)\sum\limits_{n = 1}^\infty  {\frac{{{\pi ^2}}}
{{{n^{2k}}{{\sinh }^2}\left( {\pi n} \right)}}}  + 2k{\left( { - 1} \right)^k}\sum\limits_{n = 1}^\infty  {\frac{{\pi \coth \left( {\pi n} \right)}}
{{{n^{2k + 1}}}}}  \nonumber \\
 &  + \zeta \left( {2k + 2} \right)\left\{ {2\left( {1 - {2^{ - 2k - 1}}} \right){{\left( { - 1} \right)}^{k + 1}} - 1} \right\} \nonumber \\
  & + 2\sum\limits_{\scriptstyle {k_1} + {k_2} = k,\hfill \atop \scriptstyle {k_1}, {k_2} \geqslant 1\hfill} {{{\left( { - 1} \right)}^{{k_2}}}\left( {\left( {1 - {2^{1 - 2{k_1}}}} \right){{\left( { - 1} \right)}^{{k_1}}} - 2} \right)\zeta \left( {2{k_1}} \right)\bar \zeta \left( {2{k_2}} \right)}  \nonumber \\
 &  - 4\sum\limits_{\scriptstyle {k_1} + {k_2}+k_3 = k,\hfill \atop \scriptstyle {k_1}, {k_2},k_3 \geqslant 1\hfill} {{{\left( { - 1} \right)}^{{k_1} + {k_2}}}\bar \zeta \left( {2{k_1}} \right)\bar \zeta \left( {2{k_2}} \right)\zeta \left( {2{k_3}} \right)}  = 0,
\end{align}
\begin{align}\label{4.4}
&\sum\limits_{n = 1}^\infty  {\frac{{{\pi ^2}{{\coth }^2}\left( {\pi n} \right)}}
{{{n^{2k}}}}}  + {\left( { - 1} \right)^k}\sum\limits_{n = 1}^\infty  {\frac{{{\pi ^2}}}
{{{n^{2k}}{{\sinh }^2}\left( {\pi n} \right)}}}  + 2k{\left( { - 1} \right)^k}\sum\limits_{n = 1}^\infty  {\frac{{\pi \coth \left( {\pi n} \right)}}
{{{n^{2k + 1}}}}}\nonumber \\
& + \zeta \left( {2k + 2} \right)\left( {2{{\left( { - 1} \right)}^k} - 1} \right) + 2\sum\limits_{{k_1} + {k_2} = k + 1} {\left( {{{\left( { - 1} \right)}^{{k_1} + {k_2}}} + 2{{\left( { - 1} \right)}^{{k_2}}}} \right)\zeta \left( {2{k_1}} \right)\zeta \left( {2{k_2}} \right)} \nonumber\\
& - 4\sum\limits_{{k_1} + {k_2} + {k_3} = k + 1} {{{\left( { - 1} \right)}^{{k_1} + {k_2}}}\zeta \left( {2{k_1}} \right)\zeta \left( {2{k_2}} \right)\zeta \left( {2{k_3}} \right)}  = 0.
\end{align}
Letting $k=1$ in (\ref{4.2}), (\ref{4.4}) and $k=2$ in (\ref{4.3}), we have
\begin{align}\label{4.5}
\sum\limits_{n = 1}^\infty  {\frac{{{\pi ^2}}}
{{{n^2}{{\sinh }^2}\left( {\pi n} \right)}}{{\left( { - 1} \right)}^n}}  - 2\sum\limits_{n = 1}^\infty  {\frac{\pi }
{{{n^3}\sinh \left( {\pi n} \right)}}}  - \sum\limits_{n = 1}^\infty  {\frac{{{\pi ^2}\cosh \left( {\pi n} \right)}}
{{{n^{2k}}{{\sinh }^2}\left( {\pi n} \right)}}}  + \frac{{11}}
{8}\zeta \left( 4 \right) = 0,
\end{align}
\begin{align}\label{4.6}
\sum\limits_{n = 1}^\infty  {\frac{{{\pi ^2}}}
{{{n^4}{{\sinh }^2}\left( {\pi n} \right)}}}  + 2\sum\limits_{n = 1}^\infty  {\frac{{\pi \coth \left( {\pi n} \right)}}
{{{n^5}}}}  - \frac{{13}}
{2}\zeta \left( 6 \right) = 0,
\end{align}
\begin{align}\label{4.7}
\sum\limits_{n = 1}^\infty  {\frac{{{\pi ^2}{{\coth }^2}\left( {\pi n} \right)}}
{{{n^2}}}}  - \sum\limits_{n = 1}^\infty  {\frac{{{\pi ^2}}}
{{{n^2}{{\sinh }^2}\left( {\pi n} \right)}}}  - 2\sum\limits_{n = 1}^\infty  {\frac{{\pi \coth \left( {\pi n} \right)}}
{{{n^3}}}}  - 8\zeta \left( 4 \right) = 0,
\end{align}
Using the well known identities (\cite{ES1957})
\[\sum\limits_{n = 1}^\infty  {\frac{{{\pi ^2}}}
{{{n^2}{{\sinh }^2}\left( {\pi n} \right)}}}  = \frac{2}
{3}G - \frac{{11}}
{{30}}\zeta \left( 2 \right)\quad{\rm and}\quad\sum\limits_{n = 1}^\infty  {\frac{{\pi \coth \left( {\pi n} \right)}}
{{{n^3}}}}  = \frac{7}
{2}\zeta \left( 4 \right),\]
we can get
\[\sum\limits_{n = 1}^\infty  {\frac{{{\pi ^2}{{\coth }^2}\left( {\pi n} \right)}}
{{{n^2}}}}  = \frac{{19}}
{2}\zeta \left( 4 \right) + 4\zeta \left( 2 \right)G,\]
where $G$ denotes the Catalan's constant, which is related to the complete elliptic integral
\[{\bf K} \equiv {\bf K}\left( k \right) \equiv \int\limits_0^{\frac{\pi }
{2}} {\frac{{dx}}
{{\sqrt {1 - {k^2}{{\sin }^2}x} }}} \]
by the expression
\[G = \frac{1}
{2}\int\limits_0^1 {{\bf K}dx}  = \sum\limits_{n = 1}^\infty  {\frac{{{{\left( { - 1} \right)}^n}}}
{{{{\left( {2n + 1} \right)}^2}}} = 0.915965 \ldots .} \]

\section{Two results involving harmonic number and hyperbolic function}

In this section, we give two formulae of infinite series involving harmonic number and hyperbolic function. Here the harmonic number $H_n$ is defined by \cite{B1985}
\[{H_n} := \sum\limits_{k = 1}^n {\frac{1}{k}}.\]
The generalized harmonic numbers $H^{(p)}_n$ are defined by \cite{FS1998}
\[H_n^{\left( p \right)} = \sum\limits_{k = 1}^n {\frac{1}
{{{k^p}}}} ,\quad p\in\N.\]
\begin{thm}\label{thm5.1} For positive integer $p$, then the following identity holds:
\begin{align}
&\left( {1 + {{\left( { - 1} \right)}^{p + 1}}} \right)\sum\limits_{n = 1}^\infty  {\frac{{\pi {H_n}}}{{{n^p}\sinh \left( {n\pi } \right)}}}  - \left( {p + {{\left( { - 1} \right)}^{p + 1}}} \right)\sum\limits_{n = 1}^\infty  {\frac{\pi }{{{n^{p + 1}}\sinh \left( {n\pi } \right)}}}  - \sum\limits_{n = 1}^\infty  {\frac{{{\pi ^2}\cosh \left( {n\pi } \right)}}{{{n^p}{{\sinh }^2}\left( {n\pi } \right)}}} \nonumber\\
&\quad + \sum\limits_{n = 1}^\infty  {{{\left( { - 1} \right)}^n}\frac{{\pi \coth \left( {n\pi } \right)}}{{{n^p}}}{i^{p + 1}}\left( {\psi \left( {ni} \right) + {{\left( { - 1} \right)}^{p + 1}}\psi \left( { - ni} \right) + \left( {1 + {{\left( { - 1} \right)}^{p + 1}}} \right)\gamma } \right)}\nonumber \\
&\quad + \zeta \left( {p + 2} \right)\left\{ {\left( {1 + {{\left( { - 1} \right)}^p}} \right)\left( {1 - {2^{ - p - 1}}} \right){{\left( { - 1} \right)}^{\left[ {\frac{p}{2}} \right] + 1}} - 2 - {{\left( { - 1} \right)}^p}} \right\}\nonumber\\
&\quad - 2\left( {1 + {{\left( { - 1} \right)}^p}} \right)\sum\limits_{\scriptstyle 2{n_1} + 2{n_2} = p + 2, \hfill \atop
  \scriptstyle {n_1},{n_2} \ge 1 \hfill} {\left( {1 - {2^{1 - 2{n_1}}}} \right){{\left( { - 1} \right)}^{{n_1}}}\zeta \left( {2{n_1}} \right)\zeta \left( {2{n_2}} \right)}\nonumber \\
&\quad + 2\sum\limits_{\scriptstyle 2{n_1} + {n_2} = p + 1, \hfill \atop
  \scriptstyle {n_1},{n_2} \ge 1 \hfill} {\left( {1 - \left( {1 - {2^{1 - 2{n_1}}}} \right){{\left( { - 1} \right)}^{{n_1}}}} \right)\zeta \left( {2{n_1}} \right)\zeta \left( {{n_2} + 1} \right)} \nonumber\\
&\quad + 4\sum\limits_{\scriptstyle 2{n_1} + 2{n_2} + {n_3} = p+1, \atop
  \scriptstyle n{_1},{n_2},{n_3} \ge 1} {\left( {1 - {2^{1 - 2{n_1}}}} \right){{\left( { - 1} \right)}^{{n_1}}}\zeta \left( {2{n_1}} \right)\zeta \left( {2{n_2}} \right)\zeta \left( {{n_3} + 1} \right)}\nonumber \\
&=0,
\end{align}
where The $\psi$ function is the logarithmic derivative of the Gamma function,
\[\psi \left( s \right) = \frac{d}{{ds}}\ln \Gamma \left( s \right) =  - \gamma  - \frac{1}{s} + \sum\limits_{n = 1}^\infty  {\left( {\frac{1}{n} - \frac{1}{{n + s}}} \right)}\]
and it satisfies the complement formula
\[\psi \left( s \right) - \psi \left( { - s} \right) =  - \frac{1}{s} - \pi \cot \left( {\pi s} \right).\]
Here $\Gamma \left( z \right) := \int\limits_0^\infty  {{e^{ - t}}{t^{z - 1}}dt} ,\;{\mathop{\Re}\nolimits} \left( z \right) > 0$ is called gamma function, and $\gamma$ denotes the Euler-Mascheroni constant defined by
\[\gamma  := \mathop {\lim }\limits_{n \to \infty } \left( {\sum\limits_{k = 1}^n {\frac{1}{k}}  - \ln n} \right) =  - \psi \left( 1 \right) \approx {\rm{ 0 }}{\rm{. 577215664901532860606512 }}....\]
\end{thm}
\pf From \cite{FS1998}, we have
\begin{align*}
&\psi \left( { - z} \right) + \gamma \mathop  = \limits^{z \to n} \frac{1}{{z - n}} + {H_n} + \sum\limits_{k = 1}^\infty  {\left( {{{\left( { - 1} \right)}^k}{H^{(k+1)}_n} - \zeta \left( {k + 1} \right)} \right){{\left( {z - n} \right)}^k}} ,\quad n \ge 0, \\
&\psi \left( { - z} \right) + \gamma \mathop  = \limits^{z \to  - n} {H_{n - 1}} + \sum\limits_{k = 1}^\infty  {\left( {{H^{(k+1)}_{n - 1}} - \zeta \left( {k + 1} \right)} \right){{\left( {z + n} \right)}^k}} ,\quad n > 0.
\end{align*}
Then applying the kernel $\pi \cot \left( {\pi z} \right){\left( {\psi \left( { - z} \right) + \gamma } \right)}$ to the base function
\[\frac{\pi}
{{\sinh \left( {\pi z} \right){z^p}}},\quad p\in\N,\]
namely, we need to consider the function
\[H\left( z \right) = \frac{{{\pi ^2}\cot \left( {\pi z} \right)\left( {\psi \left( { - z} \right) + \gamma } \right)}}
{{\sinh \left( {\pi z} \right){z^p}}}.\]
The only singularities are poles at the integers and $\pm ni\ (n\in\N)$. Hence, a direct residue computation gives
\begin{align*}
 & {\text{Res}}\left[ {H\left( z \right),z = n} \right] =  - p\frac{\pi }
{{{n^{p + 1}}\sinh \left( {\pi n} \right)}} - \frac{{{\pi ^2}\cosh \left( {\pi n} \right)}}
{{{n^p}\sinh \left( {\pi n} \right)}} + \frac{{\pi {H_n}}}
{{{n^p}\sinh \left( {\pi n} \right)}}, \hfill \\
 & {\text{Res}}\left[ {H\left( z \right),z =  - n} \right] = {\left( { - 1} \right)^{p + 1}}\frac{{\pi {H_{n - 1}}}}
{{{n^p}\sinh \left( {\pi n} \right)}}, \hfill \\
 & {\text{Res}}\left[ {H\left( z \right),z = ni} \right] = {\left( { - 1} \right)^{n + p + 1}}\frac{{\pi \coth \left( {\pi n} \right)}}
{{{n^p}}}\left( {\psi \left( { - ni} \right) + \gamma } \right){i^{p + 1}}, \hfill \\
 & {\text{Res}}\left[ {H\left( z \right),z =  - ni} \right] = {\left( { - 1} \right)^n}\frac{{\pi \coth \left( {\pi n} \right)}}
{{{n^p}}}\left( {\psi \left( {ni} \right) + \gamma } \right){i^{p + 1}}, \hfill \\
\end{align*}
and
\begin{align*}
{\text{Res}}\left[ {H\left( z \right),z = 0} \right]=& \zeta \left( {p + 2} \right)\left\{ {\left( {1 + {{\left( { - 1} \right)}^p}} \right)\left( {1 - {2^{ - p - 1}}} \right){{\left( { - 1} \right)}^{\left[ {\frac{p}{2}} \right] + 1}} - 2 - {{\left( { - 1} \right)}^p}} \right\}\nonumber\\
&\quad - 2\left( {1 + {{\left( { - 1} \right)}^p}} \right)\sum\limits_{\scriptstyle 2{n_1} + 2{n_2} = p + 2, \hfill \atop
  \scriptstyle {n_1},{n_2} \ge 1 \hfill} {\left( {1 - {2^{1 - 2{n_1}}}} \right){{\left( { - 1} \right)}^{{n_1}}}\zeta \left( {2{n_1}} \right)\zeta \left( {2{n_2}} \right)}\nonumber \\
&\quad + 2\sum\limits_{\scriptstyle 2{n_1} + {n_2} = p + 1, \hfill \atop
  \scriptstyle {n_1},{n_2} \ge 1 \hfill} {\left( {1 - \left( {1 - {2^{1 - 2{n_1}}}} \right){{\left( { - 1} \right)}^{{n_1}}}} \right)\zeta \left( {2{n_1}} \right)\zeta \left( {{n_2} + 1} \right)} \nonumber\\
&\quad + 4\sum\limits_{\scriptstyle 2{n_1} + 2{n_2} + {n_3} = p+1, \atop
  \scriptstyle n{_1},{n_2},{n_3} \ge 1} {\left( {1 - {2^{1 - 2{n_1}}}} \right){{\left( { - 1} \right)}^{{n_1}}}\zeta \left( {2{n_1}} \right)\zeta \left( {2{n_2}} \right)\zeta \left( {{n_3} + 1} \right)}.
\end{align*}
Thus, summing these five contributions yields the statement of the theorem.\hfill$\square$\\
When $p=1$ and $p=2$, we give
\begin{align*}
&\sum\limits_{n = 1}^\infty  {\frac{{\pi {H_n}}}{{n\sinh \left( {n\pi } \right)}}}  = \frac{1}{2}\zeta \left( 3 \right) + \pi \sum\limits_{n = 1}^\infty  {\frac{1}{{{n^2}\sinh \left( {n\pi } \right)}}}  + \frac{{{\pi ^2}}}{2}\sum\limits_{n = 1}^\infty  {\frac{{\cosh \left( {n\pi } \right)}}{{n{{\sinh }^2}\left( {n\pi } \right)}}}  + \sum\limits_{n,k = 1}^\infty  {\frac{{k\pi \coth \left( {k\pi } \right)}}{{n\left( {{n^2} + {k^2}} \right)}}{{\left( { - 1} \right)}^k}} \\
&\sum\limits_{n = 1}^\infty  {\frac{\pi }{{{n^3}\sinh \left( {n\pi } \right)}}}  + \sum\limits_{n = 1}^\infty  {\frac{{{\pi ^2}\cosh \left( {n\pi } \right)}}{{{n^2}{{\sinh }^2}\left( {n\pi } \right)}}}  - \sum\limits_{n = 1}^\infty  {\frac{{\pi \coth \left( {n\pi } \right)}}{{{n^3}}}{{\left( { - 1} \right)}^n}}  - \sum\limits_{n = 1}^\infty  {\frac{{{\pi ^2}{{\coth }^2}\left( {n\pi } \right)}}{{{n^2}}}{{\left( { - 1} \right)}^n}}  = \frac{{45}}{4}\zeta \left( 4 \right).
\end{align*}

\begin{thm}\label{thm5.2} For positive integer $p$, then the following identity holds:
\begin{align}\label{5.2}
 & \left( {1 + {{\left( { - 1} \right)}^{p + 1}}} \right)\sum\limits_{n = 1}^\infty  {\frac{{\pi {H_n}}}
{{{n^p}\sinh \left( {n\pi } \right)}}{{\left( { - 1} \right)}^n}}  - \left( {p + {{\left( { - 1} \right)}^{p + 1}}} \right)\sum\limits_{n = 1}^\infty  {\frac{\pi }
{{{n^{p + 1}}\sinh \left( {n\pi } \right)}}{{\left( { - 1} \right)}^n}} \nonumber \\
  & - \sum\limits_{n = 1}^\infty  {\frac{{{\pi ^2}\cosh \left( {n\pi } \right)}}
{{{n^p}{{\sinh }^2}\left( {n\pi } \right)}}{{\left( { - 1} \right)}^n}}   + \sum\limits_{n = 1}^\infty  {\frac{{{{\left( { - 1} \right)}^n}\pi }}
{{{n^p}\sinh \left( {n\pi } \right)}}{i^{p + 1}}\left( {\psi \left( {ni} \right) + {{\left( { - 1} \right)}^{p + 1}}\psi \left( { - ni} \right) + \left( {1 + {{\left( { - 1} \right)}^{p + 1}}} \right)\gamma } \right)}  \nonumber \\
 &  + \zeta \left( {p + 2} \right)\left\{ {\left( {1 - {2^{ - p - 1}}} \right)\left( {1 + {{\left( { - 1} \right)}^p}} \right)\left( {1 + {{\left( { - 1} \right)}^{\left[ {\frac{p}
{2}} \right] + 1}}} \right) - 1} \right\} \nonumber \\
 &  - 2\sum\limits_{\scriptstyle 2{n_1} + {n_2} = p + 1, \hfill \atop
  \scriptstyle {n_1},{n_2} \ge 1 \hfill} {\left( {1 + {{\left( { - 1} \right)}^{{n_1}}}} \right)\bar \zeta \left( {2{n_1}} \right)\zeta \left( {{n_2} + 1} \right)}\nonumber \\
  & + 2\left( {1 + {{\left( { - 1} \right)}^p}} \right)\sum\limits_{\scriptstyle 2{n_1} + 2{n_2} = p + 2, \hfill \atop
  \scriptstyle {n_1},{n_2} \ge 1 \hfill}  {{{\left( { - 1} \right)}^{{n_2}}}\bar \zeta \left( {2{n_1}} \right)\bar \zeta \left( {2{n_2}} \right)}  \nonumber\\
 &  - 4\sum\limits_{\scriptstyle 2{n_1} + 2{n_2} +2n_3= p + 1, \hfill \atop
  \scriptstyle {n_1},{n_2},n_3 \ge 1 \hfill} {{{\left( { - 1} \right)}^{{n_2}}}\bar \zeta \left( {2{n_1}} \right)\bar \zeta \left( {2{n_2}} \right)\zeta \left( {{n_3} + 1} \right)}  = 0 .
\end{align}
\end{thm}
\pf The proof is based on the kernel
\[\frac{\pi{ \left( {\psi \left( { - z} \right) + \gamma } \right)} }
{{\sin \left( {\pi z} \right)}}\]
and usual residue computation, which is applied to
\[\frac{{\pi }}
{{\sinh \left( {\pi z} \right){z^p}}}.\]
By a similar argument as in the proof of Theorem \ref{thm5.1}, we can prove the Theorem \ref{thm5.2}. \hfill$\square$\\
Letting $p=1$ in (\ref{5.2}) we obtain
\begin{align*}
\sum\limits_{n = 1}^\infty  {\frac{{\pi {H_n}}}
{{{n}\sinh \left( {n\pi } \right)}}{{\left( { - 1} \right)}^n}}  =& \frac{1}
{2}\zeta \left( 3 \right) + \sum\limits_{n = 1}^\infty  {\frac{\pi }
{{{n^{2}}\sinh \left( {n\pi } \right)}}{{\left( { - 1} \right)}^n}}  + \frac{1}
{2}\sum\limits_{n = 1}^\infty  {\frac{{{\pi ^2}\cosh \left( {n\pi } \right)}}
{{{n}{{\sinh }^2}\left( {n\pi } \right)}}{{\left( { - 1} \right)}^n}}\\&  + \sum\limits_{n,k = 1}^\infty  {\frac{{{{\left( { - 1} \right)}^n}\pi n}}
{{k\left( {{k^2} + {n^2}} \right)\sinh \left( {n\pi } \right)}}} .
\end{align*}
By using the methods and techniques of the present paper, it is possible that some of similar identities involving harmonic numbers and hyperbolic functions can be proved.

\section{Other some results involving hyperbolic function}
In \cite{Apo1990}, Apostol gave the following result of infinite product  
\begin{align}\label{6.1}
\frac{{P\left( x \right)}}
{{P\left( x^{-1} \right)}} = {e^{\frac{\pi }
{{12}}\left( {x - \frac{1}
{x}} \right)}}{x^{ - \frac{1}
{2}}},
\end{align}
where the product $P(x)$ is defined by
\[P\left( x \right): = \prod\limits_{n = 1}^\infty  {\left( {1 - {e^{ - 2\pi nx}}} \right)} .\]
Hence, from (\ref{6.1}) we obtain the closed form of infinite products
\begin{align}
&  \prod\limits_{n = 1}^\infty  {\left( {1 + {e^{ - \sqrt 2 \pi n}}} \right)}  = {2^{ - \frac{1}
{4}}}{e^{\frac{\pi }
{{12\sqrt 2 }}}}, \hfill \\
 & \prod\limits_{n = 1}^\infty  {\left( {1 + {e^{ - \pi \left( {2n - 1} \right)}}} \right)}  = {2^{\frac{1}
{4}}}{e^{-\frac{\pi }
{{24}}}}, \hfill \\
 & \prod\limits_{n = 1}^\infty  {\left( {1 + {e^{ - \pi n}}} \right)}  = {2^{ - \frac{1}
{8}}}{e^{\frac{\pi }
{{24}}}}, \hfill \\
 & \prod\limits_{n = 1}^\infty  {\left( {1 + {e^{ - 2\pi n}}} \right)}  = {2^{ - \frac{3}
{8}}}{e^{\frac{\pi }
{{12}}}}. 
\end{align}
In \cite{B1989}, Berndt proved that
\begin{align}
\sum\limits_{n = 1}^\infty  {\frac{1}
{{n\left( {{e^{2\pi n}} - 1} \right)}}}  = \frac{1}
{4}\log \left( {\frac{4}
{\pi }} \right) - \frac{\pi }
{{12}} + \log \Gamma \left( {\frac{3}
{4}} \right),
\end{align}
which is equivalent to 
\begin{align}
\prod\limits_{n = 1}^\infty  {\left( {1 - {e^{ - 2\pi n}}} \right)}  = \frac{{{\pi ^{\frac{1}
{4}}}{e^{\frac{\pi }
{{12}}}}}}
{{\sqrt 2 \Gamma \left( {\frac{3}
{4}} \right)}}.
\end{align}
Then applying the formula (\ref{6.1}), then we obtain
\begin{align}
&  \prod\limits_{n = 1}^\infty  {\left( {1 - {e^{ - \pi n}}} \right)}  = \frac{{{\pi ^{\frac{1}
{4}}}{e^{\frac{\pi }
{{24}}}}}}
{{{2^{\frac{3}
{8}}}\Gamma \left( {\frac{3}
{4}} \right)}}, \hfill \\
 & \prod\limits_{n = 1}^\infty  {\left( {1 - {e^{ - 4\pi n}}} \right)}  = \frac{{{\pi ^{\frac{1}
{4}}}{e^{\frac{\pi }
{6}}}}}
{{{2^{\frac{7}
{8}}}\Gamma \left( {\frac{3}
{4}} \right)}}, \hfill \\
  &\prod\limits_{n = 1}^\infty  {\left( {1 - {e^{ - \pi (4n - 2)}}} \right)}  = {2^{\frac{3}
{8}}}{e^{ - \frac{\pi }
{{12}}}}, \hfill \\
  &\prod\limits_{n = 1}^\infty  {\left( {1 - {e^{ - \pi (2n - 1)}}} \right)}  = {2^{\frac{1}
{8}}}{e^{ - \frac{\pi }
{{24}}}}.
\end{align}
Hence, by the definitions of hyperbolic cotangent and tangent, we have the relations
\begin{align}
&\sum\limits_{n = 1}^\infty  {\frac{{\coth \left( {\pi nx} \right)}}
{n}{{\left( { - 1} \right)}^{n - 1}}}  = \log 2 + 2\log \prod\limits_{n = 1}^\infty  {\left( {1 + {e^{ - 2\pi nx}}} \right)} ,\\
&\sum\limits_{n = 1}^\infty  {\frac{{\tanh \left( {\pi nx} \right)}}
{n}{{\left( { - 1} \right)}^{n - 1}}}  = \log 2 + 2\log \prod\limits_{n = 1}^\infty  {\left( {\frac{{1 + {e^{ - 4\pi nx}}}}
{{1 + {e^{ - 2\pi \left( {2n - 1} \right)x}}}}} \right)} ,
\end{align}
Thus, by combining the above related identities and using (\ref{2.11}), we deduce that
\begin{align}
 & \sum\limits_{n = 1}^\infty  {\frac{{\coth \left( {\frac{{\pi n}}
{2}} \right)}}
{n}{{\left( { - 1} \right)}^{n - 1}}}  = \frac{3}
{4}\log 2 + \frac{\pi }
{{12}}, \hfill \\
 & \sum\limits_{n = 1}^\infty  {\frac{{\coth \left( {\pi n} \right)}}
{n}{{\left( { - 1} \right)}^{n - 1}}}  = \frac{1}
{4}\log 2 + \frac{\pi }
{6}, \hfill \\
 & \sum\limits_{n = 1}^\infty  {\frac{{\coth \left( {\frac{{\pi n}}
{{\sqrt 2 }}} \right)}}
{n}{{\left( { - 1} \right)}^{n - 1}}}  = \frac{1}
{2}\log 2 + \frac{\pi }
{{6\sqrt 2 }}, \hfill \\
 & \sum\limits_{n = 1}^\infty  {\frac{{\tanh \left( {\frac{{\pi n}}
{2}} \right)}}
{n}{{\left( { - 1} \right)}^{n - 1}}}  = \frac{\pi }
{4} - \frac{1}
{4}\log 2, \hfill \\
 & \sum\limits_{n = 1}^\infty  {\frac{1}
{{n\sinh \left( {\pi n} \right)}}}  = \frac{\pi }
{{12}} - \frac{1}
{4}\log 2, \hfill \\
  &\sum\limits_{n = 1}^\infty  {\frac{1}
{{n\sinh \left( {2\pi n} \right)}}}  = \frac{\pi }
{6} - \frac{3}
{4}\log 2, \hfill \\
  &\sum\limits_{n = 1}^\infty  {\frac{1}
{{n\sinh \left( {\sqrt 2 \pi n} \right)}}}  = \frac{\pi }
{{6\sqrt 2 }} - \frac{1}
{2}\log 2, \hfill \\
 & \sum\limits_{n = 1}^\infty  {\frac{1}
{{\left( {2n - 1} \right)\sinh \left( {\left( {2n - 1} \right)\pi } \right)}}}  = \frac{1}
{8}\log 2.
\end{align}
{\bf Acknowledgments.} The author would like to thank Professor B.C. Berndt for his kind suggestions. 

 \end{document}